\documentclass[11pt]{amsart}
\usepackage{amssymb}
\usepackage{verbatim}

\begin{document}

\newtheorem{theorem}{Theorem}    
\newtheorem{proposition}[theorem]{Proposition}
\newtheorem{conjecture}[theorem]{Conjecture}
\def\theconjecture{\unskip}
\newtheorem{corollary}[theorem]{Corollary}
\newtheorem{lemma}[theorem]{Lemma}
\newtheorem{observation}[theorem]{Observation}
\theoremstyle{definition}
\newtheorem{definition}{Definition}
\newtheorem{remark}{Remark}
\newtheorem{question}{Question}
\def\thequestion{\unskip}
\newtheorem{example}{Example}
\def\theexample{\unskip}
\newtheorem{problem}{Problem}

\numberwithin{theorem}{section}
\numberwithin{definition}{section}
\numberwithin{equation}{section}
\numberwithin{remark}{section}

\def\reals{{\mathbb R}}
\def\torus{{\mathbb T}}
\def\integers{{\mathbb Z}}
\def\complex{{\mathbb C}\/}
\def\naturals{{\mathbb N}\/}
\def\distance{\operatorname{distance}\,}
\def\degree{\operatorname{degree}\,}
\def\kernel{\operatorname{kernel}\,}
\def\dim{\operatorname{dimension}\,}
\def\Span{\operatorname{span}\,}
\def\ZZ{ {\mathbb Z} }
\def\e{\varepsilon}
\def\p{\partial}
\def\rp{{ ^{-1} }}
\def\Re{\operatorname{Re\,} }
\def\Im{\operatorname{Im\,} }
\def\ov{\overline}
\def\eps{\varepsilon}
\def\lt{L^2}
\def\ltwos{\ell^2_s}

\def\scriptx{{\mathcal X}}
\def\scriptj{{\mathcal J}}
\def\scriptr{{\mathcal R}}
\def\scripts{{\mathcal S}}
\def\scriptb{{\mathcal B}}
\def\scripta{{\mathcal A}}
\def\scriptk{{\mathcal K}}
\def\scriptd{{\mathcal D}}
\def\scriptp{{\mathcal P}}
\def\scriptl{{\mathcal L}}
\def\scriptv{{\mathcal V}}
\def\scripti{{\mathcal I}}
\def\scripth{{\mathcal H}}
\def\scriptm{{\mathcal M}}
\def\scripte{{\mathcal E}}
\def\scriptt{{\mathcal T}}
\def\scriptb{{\mathcal B}}
\def\scriptf{{\mathcal F}}
\def\scriptn{{\mathcal N}}

\def\frakm{{\mathfrak m}}

\def\ndiag{{\mathcal N}_{\rm diag}}
\def\nmain{{\mathcal N}_{\rm main}}
\def\boldn{{\mathbf N}}
\def\zeroone{[0,1]}

\author{Michael Christ}
\address{
        Michael Christ\\
        Department of Mathematics\\
        University of California \\
        Berkeley, CA 94720-3840, USA}
\urladdr{math.berkeley.edu/$\sim$mchrist}
\email{mchrist@math.berkeley.edu}
\thanks{This material is based upon work supported by the 
National Science Foundation under Grant No.\ 0401260.}

\date{March 17, 2005}

\title
[Nonuniqueness of weak solutions] 
{Nonuniqueness of Weak Solutions\\ of the Nonlinear Schr\"odinger Equation}

\begin{abstract}
Generalized solutions of the Cauchy problem
for the one-dimensional periodic nonlinear Schr\"odinger equation,
with cubic or quadratic nonlinearities, 
are not unique. For any $s<0$ there exist nonzero generalized solutions
varying continuously in the Sobolev space $H^s$, 
with identically vanishing initial data.
\end{abstract}


\maketitle

\section{Introduction}

The Cauchy problem for the one-dimensional periodic cubic
nonlinear Schr\"odinger equation is
\begin{equation} \label{nlsivp}
\left\{
\begin{aligned}
&iu_t + u_{xx} + \omega |u|^2 u=0
\\
&u(0,x)=u_0(x)
\end{aligned}
\right.
\tag{NLS}
\end{equation}
where $x\in\torus=\reals/2\pi\integers$,
$t\in\reals$, and the parameter $\omega$ equals $\pm 1$.
Bourgain \cite{bourgainperiodic} has shown this problem to be
wellposed in the Sobolev space $H^s$
for all $s\ge 0$.
That is, there exists a Banach space $Y\subset C^0([0,T], H^s(\torus))
\cap L^3([0,T]\times\torus)$
such that for any $u_0\in H^s$ there exists a solution $u\in Y$,
and solutions within the class $Y$ are unique. Here $T$ depends on 
the $H^s$ norm of the initial datum.
An alternative proof of existence of solutions in $C^0([0,T], H^s(\torus))$
for $s\ge 0$,
without any uniqueness assertion, was recently given \cite{christpowerseries}.

On the other hand, the wellposedness theory breaks down in 
Sobolev spaces of negative order.
For $s<0$ the mapping from smooth data to solutions
fails to be uniformly continuous \cite{burqgerardtzvetkov}
with respect to the $H^s$ norm,
and is unstable in stronger senses \cite{cct3} as well.
For $s<-\tfrac12$, for any $\eps>0$ there exists 
a solution\footnote{The construction of these solutions
in \cite{cct3} does {\em not\/} permit any passage to the limit to 
obtain nonvanishing solutions with vanishing initial data.}
$u\in C^\infty$ satisfying $\|u\|_{C^0([0,\eps],H^s)}>\eps^{-1}$
with initial datum satisfying $\|u(0,\cdot)\|_{H^s}\le\eps$.

There remains the question of unconditional uniqueness,
that is, uniqueness of solutions belonging to 
$C^0([0,T],H^s)$, without further restrictions.
As it stands, this question is not well formulated,
because of the lack of any well-defined product for general
sufficiently singular distributions.
In particular,
the information $u\in C^0([0,T],H^s)$ alone is insufficient to 
ensure that the nonlinear expression $|u|^2 u$ 
has a natural interpretation as a space-time distribution.
When $s$ is sufficiently large this expression makes sense, and 
solutions in $C^0([0,T],H^s)$ are then well known to be unique. 
More refined work has established sufficient conditions on $s$
for unconditional uniqueness for various equations; see for instance
\cite{planchonetal} and references cited there.

In this paper we establish nonuniqueness of solutions to the Cauchy
problem for the (periodic, cubic) nonlinear Schr\"odinger equation
and its variants with quadratic nonlinearities
in classes $C^0([0,T],H^s(\torus))$ for $s<0$.
While the paper focuses on one prototypical equation
and some of its variants, the underlying construction is quite general. 
Two caveats must be admitted: (i) The solutions constructed
are sufficiently singular that the meaning of the nonlinear terms 
in the equations must be clarified before it can be discussed whether
the differential equation is actually satisfied.
We prove that the required nonlinear expressions have reasonable and canonical
interpretations, and that the differential equations hold under these
interpretations.
(ii) In the cubic case, the differential equation is modified slightly. 
The resulting modified Cauchy problem \eqref{modnlsivp} has a reasonable
existence theory with uniformly continuous dependence on initial data, 
in a natural but weak sense, 
for a {\em wider} class of function spaces than does \eqref{nlsivp}. 
See below for more precise discussions of these two points.

For \eqref{modnlsivp} there exist certain function
spaces $\scripth$ such that rather canonical solutions 
in $C^0([0,T],\scripth)$
exist for all initial data $u_0\in \scripth$,  
with uniformly continuous dependence upon
initial data, yet solutions in $C^0([0,T],\scripth)$ fail to be unique.
The same holds for the nonlinear Schr\"odinger equation with 
certain quadratic nonlinearities, in Sobolev spaces $H^s$ for all strictly
negative $s$.


\section{Results}
\subsection{Definitions}
Our modified Cauchy problem is
\begin{equation} \label{modnlsivp}
\left\{
\begin{aligned}
&iu_t + u_{xx} + \omega 
\boldn(u)=0
\\
&u(0,x)=u_0(x)
\end{aligned}
\right.
\tag{NLS$^*$}
\end{equation}
where
\begin{align}
&\boldn(u) = \big( |u|^2-2\mu(|u|^2)) u
\\
&\mu(f) = (2\pi)\rp\int_{\torus} f(x)\,dx.
\end{align}
$\mu(|u(t,\cdot)|^2)$
is independent of $t$ for all sufficiently smooth solutions;
modifying the equation in this way merely introduces a unimodular scalar
factor $e^{2i\mu t}$, where $\mu =\mu(|u_0|^2)$.
It is always assumed that $\omega\ne 0$,
so that the equation is genuinely nonlinear.
For parameters $s<0$,
$\mu(|u_0|^2)$ is not defined for typical $u_0\in H^s$, but
of course the same goes for $|u_0(x)|^2$. 
Subtracting $2\mu(|u|^2)u$ makes the
equation better behaved, as discussed below; it contributes to the nonuniqueness
of solutions by making it possible to reasonably interpret
the modified differential equation for a wider class of distributions 
than the unmodified equation, but does not directly produce any wild behavior.

We will work with the partial Fourier transform, which is defined for
smooth functions $f(t,x)$ by
\begin{equation}
\widehat{f}(t,n) = (2\pi)\rp\int_\torus f(t,x)e^{-inx}\,dx
\ \text{ for } n\in\integers,
\end{equation}
and is extended to distributions by continuity.

\begin{definition} \label{defn:Fouriercutoff}
A sequence of Fourier cutoff operators
is any sequence of linear operators
$(\scriptp_N)_{N\in\naturals}$
which act on $\scriptd'(\torus)$, and
are of the Fourier multiplier form 
$\widehat{\scriptp_N f}(n) = m_N(n)\widehat{f}(n)$ 
where
the functions $m_N:\integers\to\complex$
each have finite support,
are uniformly bounded, and 
satisfy $\lim_{N\to\infty}m_N(k)=1$
for all $k\in\integers$.
\end{definition}

Let $\scriptn$ be some nonlinear functional acting
on functions of $(t,x)$.
\begin{definition} \label{defn:truncation}
Let $u\in\scriptd'((0,1)\times\torus)$ be a distribution.
$\scriptn(u)$ is said to exist and to equal 
$v\in\scriptd'((0,1)\times\torus)$
if for every sequence $(\scriptp_N)$ of Fourier cutoff operators,
\begin{equation} \label{truncationlimit}
\lim_{N\to\infty} \scriptn(\scriptp_N u)=v
\end{equation} 
in the topology of $\scriptd'((0,1)\times\torus)$.
\end{definition}
We emphasize that \eqref{truncationlimit} is to hold for every
sequence $(\scriptp_N)$, not merely for one sequence.
Under general theories of multiplication of distributions
\cite{biagioni},\cite{colombeau}, products of the objects
discussed are always defined, but these products depend on
the choice of approximating truncation operators.
One could require still more of $u$ by replacing Fourier
cutoff operators by an appropriate class of pseudodifferential
operators implementing cutoffs in phase space rather than
merely in frequency space; we have not systematically investigated
this more restrictive notion of existence for the solutions
constructed in this paper.


\begin{definition} \label{eqnsense} \label{defn:weaksolution}
$u\in C^0(\zeroone,H^s(\torus))$
will be said to be a weak solution
of \eqref{modnlsivp} in the extended sense if
$u(0,\cdot)=u_0$, 
$\boldn(u)$ exists in the sense of Definition~\ref{defn:truncation},
and $u$ satisfies  $iu_t +u_{xx}+\boldn(u)=0$
in the distribution sense in $(0,1)\times\torus$
with this interpretation of $\boldn(u)$.
\end{definition}

See \cite{dix} for some discussion of this and related, less restrictive,
notions of weak solutions.

For any function space $\scripth=\scripth(\torus)$,
$C^{-1}(\zeroone,\scripth)$
will denote the space of all space-time distributions $F(s,x)$
such that $\tilde F(t,x) = \int_0^t F(s,x)\,ds$ belongs
to $C^0(\zeroone,\scripth)$, and $\|F\|_{C^{-1}(\zeroone,\scripth)}$
is the ``norm'' $\max_{t\in\zeroone}\|\tilde F(t,\cdot)\|_{\scripth}$.

The construction will rely on solutions of the inhomogeneous problem
\begin{equation} \label{drivenmodnlsivp}
\left\{
\begin{aligned}
&iv_t + v_{xx} + \omega \boldn(v)=F
\\
&v(0,x)=v_0(x).
\end{aligned}
\right.
\end{equation}
We refer to $F$ as a driving force.
Constructions of Scheffer \cite{scheffer} and Shnirelman \cite{shnirelman} 
of nonunique solutions for the Euler equation
have also utilized solutions of inhomogeneous equations.

\subsection{Nonuniqueness for the cubic nonlinearity}
\begin{theorem} \label{thm:first}
For any $s<0$ and $\omega\ne 0$,
there exists a space-time distribution $u\in C^0([0,1],H^s)$, 
not identically vanishing,
which is a weak solution of \eqref{modnlsivp} in the extended sense,
with initial datum $u_0\equiv 0$.
Moreover, the limit \eqref{truncationlimit} defining 
$e^{-it\Delta}\boldn(u)$ exists in the $C^{-1}(\zeroone,H^s)$ norm.
\end{theorem}
It can be shown by an elaboration of the proof that
for any initial datum with $\widehat{u_0}\in\ell^1$,
there exist $T>0$ and two distinct weak solutions in $C^0([0,T],H^s)$
of \eqref{modnlsivp}.
Similar extensions are possible for all theorems stated below.

The solution $u$ qualifies as a solution in a second sense:
There exist sequences of functions
$f_n\in C^\infty([0,1]\times\torus)$,
such that $e^{-i\Delta t}f_n(t,x)\to 0$
in $C^{-1}([0,1],H^s)$ norm as $n\to\infty$,
and solutions $u_n\in C^0(\zeroone,H^1)$
of \eqref{drivenmodnlsivp} with 
driving forces $f_n$ and initial data $u_0\equiv 0$,
such that $u_n\to u$ in $C^0(\zeroone,H^s)$ norm as $n\to\infty$.

While Theorem~\ref{thm:first} concerns rather irregular weak solutions,
the essence of the construction is the following approximation result
for smooth solutions of the inhomogeneous problem.
\begin{proposition} \label{prop:surgery}
Let $s<0$ and $\omega\ne 0$.
Suppose that $u\in C^\infty(\zeroone\times\torus)$,
and that each Fourier coefficient $\widehat{u}(t,n)$
vanishes to infinite order as $t\to 0^+$.
Then for any $\eps>0$ there exist $v,F\in C^\infty(\zeroone\times\torus)$,
each of whose Fourier coefficients vanishes to infinite order as $t\to 0^+$,
such that $v$ is a solution of the inhomogeneous Cauchy problem
\eqref{drivenmodnlsivp} with driving force $F$,
with bounds
\begin{align}
& \|v-u\|_{C^0(\zeroone,H^s)}\le\eps
\\
& \|e^{-it\Delta}F\|_{C^{-1}(\zeroone,H^s)}\le\eps.
\end{align}
\end{proposition}
The other theorems stated below are based on analogous facts.

\subsection{Earlier nonuniqueness results}
Theorem~\ref{thm:first} should be contrasted with the examples of 
Scheffer \cite{scheffer} and Shnirelman \cite{shnirelman} 
of nonunique weak solutions of the
(periodic, two-dimensional) incompressible Euler equation in $C^0([0,T],H^0)$.
The notion of a weak solution is less problematic in that framework, for
the nonlinear term $v\cdot\nabla v$ is well-defined as a space-time
distribution, under the usual straightforward 
interpretation via integration by parts, for any $v\in C^0([0,T],H^0)$.

A result related to nonuniqueness for the nonlinear Schr\"odinger
equation on the real line
has been established by Kenig, Ponce, and Vega \cite{kenigetal}:
With a Dirac mass as initial datum, either there exists no solution,
or there exists more than one solution.\footnote{\cite{kenigetal} does not
address the issue of defining $|u|^2 u$, and the number of solutions 
could conceivably depend on the definition used. What is actually proved is that
for any interpretation of $|u|^2 u$ that is 
appropriately invariant under Galilean symmetries of the equation, 
there exists either no solution, or more than one solution.}
Dix \cite{dix} has shown nonuniqueness of weak solutions in 
$C^0(H^s)$ for Burgers' equation, for $s<-\tfrac12$,
via the Cole-Hopf transformation,
which transforms solutions of the heat equation to solutions of Burgers' equation
by taking a logarithm. 

\subsection{Nonuniqueness in more restrictive function spaces}
We will also establish, by a slightly more complicated argument, 
the analogue of Theorem~\ref{thm:first} for certain less standard function spaces.
These are the spaces $\scripth^{p}$
for $p\in [1,\infty)$,  defined by
\begin{definition}
$\scripth^{p}(\torus)=\{f\in \scriptd(\torus): 
\widehat{f}(\cdot)\in\ell^p\}$.
\end{definition}
\noindent
Here $\scriptd(\torus)$ is the usual space of distributions, and $\scripth^{p}$
is equipped with the norm
$\|\widehat{f}\|_{\ell^{p}(\integers)}$.

The Cauchy problem \eqref{modnlsivp}
in $\scripth^p$ exhibits certain attributes of wellposedness
for all $p\in[1,\infty)$ 
\cite{christpowerseries}:
For any $R<\infty$ there exists $T>0$ such that
the solution operator $u_0\mapsto u(t,x)$,
defined initially for all $u_0\in H^1$,
is uniformly continuous (even real analytic) as a mapping
from $\{u_0\in H^1: \|u_0\|_{\scripth^p}\le R\}$,
equipped with the  $\scripth^p$ topology,
to $C^0([0,T],\scripth^p)$. Moreover the mapping $u_0\mapsto u$
defined by extending this mapping from the dense subspace 
to all of $\scripth^p$
is actually real analytic, and the function $u(t,x)$ thus defined
is a weak solution of the differential equation in 
the extended sense.
The unmodified Cauchy problem \eqref{nlsivp} lacks these features
for all $p>2$; the modified equation is better behaved.

\begin{theorem} \label{thm:second}
Let $p>2$ and $\omega\ne 0$.
There exists a weak solution $u\in C^0([0,1],\scripth^p)$ of \eqref{modnlsivp},
in the extended sense,
which does not vanish identically but has initial datum $u_0\equiv 0$.
Moreover, the limit \eqref{truncationlimit} defining 
$e^{-it\Delta}\boldn(u)$ exists in the $C^{-1}(\zeroone,\scripth^p)$ norm.
\end{theorem}

\subsection{Quadratic nonlinearities}
Consider next the Cauchy problem 
\begin{equation} \label{quadivp}
\left\{
\begin{aligned}
&iu_t + u_{xx} + \omega 
Q(u)=0
\\
&u(0,x)=u_0(x)
\end{aligned}
\right.
\tag{NLS$_2$}
\end{equation}
where 
\begin{equation} \label{Qs}
Q(u)=u^2,\  = \bar u^2,\  \text{or } =|u|^2-\mu(|u|^2).
\end{equation}

\begin{theorem}
Let $s<0$ and $\omega\ne 0$.
For the Cauchy problem \eqref{quadivp} with 
any of the nonlinearities\footnote{
The same conclusion holds for cubic nonlinearities
$u^3$ and $\bar u^3$; no modification of the nonlinearity
like that in \eqref{modnlsivp} is required.}
\eqref{Qs}, there exists $u\in C^0([0,1],H^s)$ 
which is a weak solution in the extended sense, 
does not vanish identically, and has
initial datum $u_0\equiv 0$. 
Moreover, $\lim_{N\to\infty} e^{-it\Delta}Q(\scriptp_N u)$ 
exists in the $C^{-1}(\zeroone,H^s)$ norm
for any sequence of operators $\scriptp_N$
satisfying the conditions of Definition~\ref{defn:Fouriercutoff}.
\end{theorem}
\noindent

For $Q=u^2$ or $\bar u^2$, this Cauchy problem is wellposed in 
$H^s$ for all $s>-\tfrac12$ \cite{kpvquad}, in the usual sense; 
for any initial datum in $H^s$ there exists a solution belonging
to a space more restrictive than $C^0(\zeroone,H^s)$, and
within this smaller space the solution is unique.
Thus for $s\in (-\tfrac12,0)$ we have simultaneously wellposedness in $H^s$
in the usual sense, and nonuniqueness of weak solutions in the extended sense
in $C^0(\zeroone,H^s)$.

\subsection{Discussion}
The construction proceeds as follows.
We consider a sequence of exact solutions $u_\nu$
of the modified Cauchy problem with initial data zero and
with driving forces $f_\nu = \sum_{|k|\ge M_\nu} c_{k,\nu}(t)e^{ikx}$,
where $M_\nu\to\infty$ as $\nu\to\infty$. To leading order, $f_\nu$
contributes $v_\nu(t) = -i\int_0^t e^{i(t-s)\Delta}
f_\nu(s)\,ds$ to the solution $u_\nu$. 
We choose $f_{\nu+1}$ so that $\boldn(v_{\nu+1})\approx f_\nu$,
modulo a very small remainder; it is essential to
work in function spaces $\scripth(\torus)$
in which it is possible to simultaneously make $v_{\nu+1}$
small in $C^0(\zeroone,\scripth)$, and $\boldn(v_{\nu+1})$ large
in $C^{-1}(\zeroone,\scripth)$.
Thus nonuniqueness arises via an infinite cascade
of ``energy'' from high spatial Fourier modes to lower Fourier modes,
that is, from small spatial scales to large scales.
Our construction and that of Shnirelman \cite{shnirelman}
have in common both
the use of driving forces tending weakly to zero,
and the exploitation of this reverse energy cascade.

The motivation for the construction is that if the evolution is viewed
as a coupled system of ordinary differential equations for
the spatial Fourier coefficients of $\widehat{u}(t,n)$,
then because this system has infinite dimension, 
uniqueness should be expected to fail without
some growth restriction as $|n|\to\infty$. 
The main issues in the construction are then that
exponential growth with respect to $n$ must be avoided, 
and that the inverse
energy cascade inevitably produces many undesired terms along with
terms useful in the construction, and it is required to make all undesired
terms small in order to keep the $H^s$ norm finite,
while useful terms are large and prescribed.

\subsection{Extensions, and potential extensions}
Various related results follow in a straightforward way from the same method.
\begin{itemize}
\item
Let $L$ be any linear operator of the form $\widehat{Lu}(n) 
= \sigma(n) \widehat{u}(n)$ where $\sigma$ is real-valued.
Then Theorem \ref{thm:first} and its proof
carry through, nearly verbatim, when the linear term $u_{xx}$
in the differential equation is replaced by $Lu$. 
More generally, if $\sigma$ has nonnegative imaginary part,
the construction goes through if rewritten without the substitution
\eqref{unitaryconjugation}.
\item
Generalization to higher dimensions is likewise straightforward.
\item
Many other nonlinearities can be treated by the same method.
Suitable modifications, analogous to
the subtraction of $2\mu(|u|^2)u$, are often needed
in order to make sense of the equation in $H^s$ for negative $s$.
\item
In particular, the periodic Korteweg-de Vries equation admits nonunique
solutions, in the extended weak sense, in $C^0(H^s)$ for all $s<0$. 
This contrasts with the
work of Kappeler and Topalov \cite{kappelertopalovkdv}, who
have proved existence of quite canonical solutions in $C^0(H^s)$, 
which depend continuously on initial data in $H^s$ for  all $s>-1$.
These ``solutions'' were only proved to satisfy the PDE in
the quite weak sense of being limits in $C^0(H^s)$ of $C^\infty$ solutions.
Our construction shows that if this notion of solution is 
liberalized by allowing limits of smooth solutions of inhomogeneous
Cauchy problems with smooth driving forces tending to zero in 
the natural space $e^{it\Delta}(C^{-1}(H^s))$,
then solutions are no longer unique.
\item
The construction applies to semilinear hyperbolic equations
$u_{tt}-\Delta_x u + \scriptn(u)=0$, 
for many nonlinearities $\scriptn$.
\end{itemize}

Other extensions and variants are at present more speculative.
It appears to be possible to:
\begin{itemize}
\item
Sharpen the examples of Scheffer
and Shnirelman for the Euler equation via this construction, 
to produce solutions in $C^0(L^2)$ 
rather than merely in $L^2([0,T]\times\torus)$.
\item
Establish nonuniqueness of the initial value problem for the
Navier-Stokes equation, for solutions in the extended weak sense
in $C^0(H^s)$ for $s$ strictly negative. 
This does not address the question of uniqueness of Leray's weak solutions 
in $C^0(H^0)$.
\item
Extend the construction to positive Sobolev exponents,
for a certain class of artificial equations such as
$iu_t+u_{xx} + u_x\bar u=0$. 
\end{itemize}
However, at present none of this has been verified in detail.

One feature of the
construction is that it is relatively insensitive to the degree
of the (semilinear) nonlinear term, in contrast to 
the behavior of threshold exponents in wellposedness theorems. 

I thank Betsy Stovall for proofreading the manuscript.


\section{Reformulation as an ordinary differential equation}
We reformulate the Cauchy problem \eqref{modnlsivp} as an infinite
coupled system of ordinary differential equations for the Fourier coefficients
of $u$.
Define
\begin{equation}
\sigma(j,k,l,n) = n^2-j^2+k^2-l^2.
\end{equation}
Written in terms of Fourier coefficients $\widehat{u}_n(t) = \widehat{u}(t,n)$
and $\widehat{F}_n(t) = \widehat{F}(t,n)$,
the differential equation 
$iu_t + u_{xx} + \omega \boldn u=F$
becomes
\begin{equation} \label{firstode}
i\frac{d \widehat{u}_n}{dt} -n^2 \widehat{u}_n
+ \omega \sum_{j-k+l=n} \widehat{u}_j\overline{\widehat{u}_k}\widehat{u}_l
-2\omega \widehat{u}_n \sum_m |\widehat{u}_m|^2
=\widehat{F_n}(t).
\end{equation}
Here the first summation is taken over all $(j,k,l)\in\integers^3$
satisfying the indicated identity, and the second over all $m\in\integers$.
The term $-2\omega \widehat{u}_n\sum_m |\widehat{u}_m|^2$
cancels out certain terms of the first sum. Eliminating these and
substituting\footnote{This substitution is natural but does not
materially simplify the analysis here. For dissipative equations
it should of course be avoided.}
\begin{equation} \label{unitaryconjugation}
y_n(t) = e^{in^2 t}\widehat{u}(t,n),
\end{equation}
\eqref{firstode} becomes
\begin{equation} \label{secondode}
\frac{d y_n}{dt} = i\omega \sum_{j-k+l=n}^* y_j \bar y_k y_l e^{i \sigma(j,k,l,n)t}
-i\omega |y_n|^2 y_n
-i e^{in^2 t}\widehat{F_n}(t) 
\end{equation}
where the notation $\sum_{j-k+l=n}^*$ means that the sum is taken over all
$(j,k,l)\in\integers^3$ for which neither $j=n$ nor $l=n$.

For a sequence $a$ define
\begin{equation}
\|a\|_{\ltwos}^2 = \sum_{n\in\integers} |a_n|^2(1+n^2)^s.
\end{equation}
Clearly $y\in C^0([0,T],\ltwos)$ if and only if $u\in C^0([0,T],H^s)$,
with identical norms.

For any complex-valued sequence $z$
define
$\ndiag(t)(z)$ and $\nmain(t)(z)$ 
to be the sequences whose $n$-th terms are 
\begin{align*}
&[\nmain(t)(z)]_n =
i\omega \sum_{j-k+l=n}^* z_j \bar z_k z_l e^{i \sigma(j,k,l,n)t}
\\
&[\ndiag(t)(z)]_n 
=-i\omega |z_n|^2 z_n.
\end{align*}
and define
\begin{equation*}
\scriptn(z) = \nmain(z) + \ndiag(z).
\end{equation*}
For each $t\in\reals$,
$\scriptn(t)$ is a nonlinear operator which
acts on a numerical sequence $z=(z_n)_{n\in\integers}$,
and produces another numerical sequence.

We will work with sequence-valued functions $y$ of $t$,
and $\scriptn(y)$ will denote the sequence-valued function
$\scriptn(t)(z)$ where $z=y(t)$.
With this notation, \eqref{secondode} becomes
\begin{equation} \label{integraleqn}
\frac{dy}{dt} = \scriptn(y) + f
\end{equation}
where 
\begin{equation*}
f_n(t) = -ie^{in^2 t}\widehat{F_n}(t).
\end{equation*}

We say that a sequence-valued function $h(t)=(h_n(t))_{n\in\integers}$ 
of $t\in[0,1]$ has support contained in $S\subset\integers$
if $h_n(t)\equiv 0$ for all $t\in[0,1]$, for every $n\notin S$.
Thus we may speak of sequence-valued functions with finite supports.

\section{The main step}

Expressed in terms of Fourier coefficients, Proposition~\ref{prop:surgery}
becomes
\begin{proposition} \label{prop:forcing}
Let $s<0$.
Let $x\in C^\infty([0,1])$
be a finitely supported sequence-valued function 
such that
$x(t,\cdot)$ vanishes to infinite order as $t\to 0^+$.
Then for any $\eps>0$ 
there exist finitely supported sequence-valued functions
$y,g\in C^\infty([0,1])$ satisfying
\begin{equation} \label{eqnfory}
\left\{
\begin{aligned}
& \frac{dy}{dt} = \scriptn(y) + g(t)
\\
&y(t,\cdot) \ \text{vanishes to infinite order as $t\to 0^+$}
\end{aligned}
\right.
\end{equation}
with
\begin{align}
&\|y-x\|_{C^0([0,1],\ltwos)}\le\eps
\label{boundfory}
\\
&\|g\|_{C^{-1}([0,1],\ltwos)}\le\eps.
\label{boundforg}
\end{align}
\end{proposition}
\noindent
Moreover, for any $M<\infty$, 
$y$ may be constructed so that $y-x$ and $g$ are supported in $[M,\infty)$.

\subsection{Construction of $y$}
Define 
\begin{equation}f = \frac{dx}{dt}-\scriptn(x).\end{equation}
Since $x$ has finite support, so does $f$.
Let $S$ be a finite set in which $f$ is supported,
and write $S=\{n_j: 1\le j\le A\}$ where the $n_j$ are distinct.
Choose a finite set $S^{\dagger} \subset\integers\cap [M,\infty)$, 
as follows. 
First choose $m_1\ge M$, and define $m'_1$ by the equation
$2m_1-m'_1=n_1$. Make $m_1$ sufficiently large to ensure
that $m'_1\ge M$ as well.
Choose $m_2\ge M$ very large relative to $m_1,m'_1$,
and define $m'_2$ by $2m_2-m'_2=n_2$. 
Then choose $m_3,m'_3,m_4,m'_4\dots$ in that order, satisfying 
\begin{equation}
2m_j-m'_j=n_j
\ \text{for all $1\le j\le A$,} 
\end{equation}
and let $S^{\dagger} =\{m_1,m'_1,\cdots,m_A,m'_A\}$. 
The elements of $S^\dagger$ are to be chosen to satisfy 
additional constraints:
\begin{enumerate}
\item  \label{itemcubic}
If $k,l,m\in S^\dagger$ 
and if $l\notin\{k,m\}$ then $|k-l+m|\ge M$
unless $(k,l,m)=(m_j,m'_j,m_j)$ for some $j$.
\item  \label{itemquadratic}
If $k,l\in S^\dagger$ and $n$ belongs to the support of $x$
then $|k-n+l|\ge M$. Moreover $|k-l+n|\ge M$ provided that $k\ne l$.
\item  \label{itemlinear}
If $k\in S^\dagger$ and $m,n$ belong to the support of $x$
then $|k-m+n|\ge M$ and $|m+k-n|\ge M$.
\end{enumerate}
Since each $m'_j$ is approximately twice as large as $m_j$,
and since the support of $x$ is finite,
all these conditions will hold, provided that $m_1$ is sufficiently large
and each subsequent $m_j$
is chosen sufficiently large relative to $m_1,\cdots,m_{j-1}$,
while $m'_j$ is defined to be $2m_j-n_j$.

Choose $C^\infty$ functions $\{h_m(t): m\in S^\dagger\}$ 
that vanish to infinite order as $t\to 0$ 
and satisfy
\begin{equation}
i\omega \overline{h_{m'_j}}(t)h_{m_j}^2(t)
\equiv \tfrac12 e^{-i\sigma(m_j,m'_j,m_j,n_j)t}
f_{n_j}(t)
\end{equation}
for each $n_j\in S$.
It is essential that these functions be chosen so that 
$\max_{m\in S^\dagger} \|h_m\|_{C^0(\zeroone)}$ is bounded above
by a finite quantity depending only\footnote{The $C^1$ norms
of the functions $h_m$ will be finite but must depend on $S^\dagger$.
This is due to the dispersive nature of the PDE, and prevents us from making
$dh/dt$ small in $C^0(\ltwos)$. This is an essential part of the obstruction
to extending the construction to positive Sobolev exponents.}
on $S$ and on $f$, not on the choice of $S^{\dagger}$ itself.
Define $h=(h_j(t))_{j\in\integers}$
by $h_j(t)=0$ for all $j\notin S^\dagger$,
and $h_j$ as above for all $j\in S^\dagger$.
Define
\begin{equation} y=x+h.\end{equation}

\subsection{Remainder terms}
Define 
\begin{equation}
g = \frac{dy}{dt}-\scriptn(y). 
\end{equation}
Since $x,h$ have disjoint supports,
$\ndiag(x+h)=\ndiag(x)+\ndiag(h)$. 
Consequently
\begin{equation} \label{gformula}
g = \big(f-\nmain(h)\big)
+ \frac{dh}{dt} - \ndiag(h)
- \big( \nmain(x+h)-\nmain(x)-\nmain(h) \big).
\end{equation}

The bounds on $y-x$ and $g$ in Proposition~\ref{prop:forcing}
will now be established.
As in other constructions of poorly behaved solutions 
\cite{burqgerardtzvetkov},\cite{cct1},\cite{cct2},\cite{cct3},
we work in a regime in which nonlinear effects are more powerful
than dispersion.
\begin{lemma} \label{mainerrorbound}
Let $x$ be as in the hypotheses of Proposition~\ref{prop:forcing},
and let $h$ be constructed as above.
Then for any $\eps>0$ there exists $M<\infty$ such
that if $S^\dagger$ is chosen as specified, then
\begin{align}
&\|h\|_{C^0(\zeroone,\ltwos)}\le\eps
\\
&\|\frac{dh}{dt}\|_{C^{-1}(\zeroone,\ltwos)}\le\eps  
\\
&\|\nmain(x)+\nmain(h) -\nmain(x+h)\|_{C^{0}(\zeroone,\ell^2_s)}
\le\eps
\\
&\|\nmain(h)-f\|_{C^{0}(\zeroone,\ell^2_s)}
\le\eps
\\
&\|\ndiag(h)\|_{C^{0}(\zeroone,\ell^2_s)}\le\eps
\end{align}
\end{lemma}

\begin{proof}
$h_k(t)$ vanishes for all $k\notin S^\dagger$, and is bounded uniformly
by a finite constant depending on $f$, independent of the choice
of $S^\dagger$. The cardinality of $S^\dagger$ likewise depends only
on $x$.
Since $s$ is strictly negative, it follows that
$\|h\|_{C^0(\zeroone,\ltwos)}\le CM^{s}$. 

The bound for $\frac{dh}{dt}$ is merely a restatement of the bound
for $h$.
$\ndiag(h)$ is also supported in $S^\dagger$, and the same reasoning
as for $h$ applies to it. 

$\nmain(h)-f$ is supported on $\{n: |n|\ge M\}$,
by  \ref{itemcubic}. 
The term $\nmain(x)+\nmain(h) -\nmain(x+h)$
is supported in the same set, by \ref{itemquadratic}
and \ref{itemlinear}.
Therefore the same reasoning applies to them and
yields the same bound $CM^s$. \end{proof}

\section{A solution with zero initial datum}

\subsection{Construction of the solution}
By induction on $n\in\naturals$,
we construct a sequence of finitely supported $C^\infty$
sequence-valued functions $x^{(n)}\in C^1(\zeroone)$ 
which vanish to infinite order as $t\to 0$.
To begin, choose $x^{(1)}$ to be smooth, to have finite support, 
to vanish to infinite order as $t\to 0$,
and moreover to have $0$-th component satisfying
\begin{equation} \label{initialbigness}
\|x^{(1)}_0\|_{C^0(\zeroone)}\ge 1.
\end{equation}

For the inductive step, construct 
$x^{(n+1)}=y$
by applying Proposition~\ref{prop:forcing} to $x=x^{(n)}$.
Define the increments $h^{(n)} = x^{(n+1)}-x^{(n)}$ and the
driving forces $f^{(n+1)} = \frac{dx^{(n)}}{dt}-\scriptn(x^{(n)})$.
Then by induction $h^{(n)}$ and hence $x^{(n+1)}$ vanish to infinite order
as $t\to 0$.
Taking $\eps$ to be sufficiently small in the conclusion of the proposition
at each step, we obtain bounds
\begin{align}
& \|h^{(n)}\|_{C^{0}(\zeroone,\ltwos)}\le\delta_n
\label{itiscauchy}
\\
&\|f^{(n+1)}\|_{C^{-1}(\zeroone,\ltwos)}
\le \delta_{n}
\\
&\delta_n\le 2^{-n-1}
\end{align}
and moreover each $\delta_n$ may be arranged to be as small as may be desired, 
relative to any quantity depending only on $x^{(n)}$.
Moreover $h^{(n)}$ and $f^{(n+1)}$ are naturally expressed as finite
sums of various constituent quantities, discussed in the proof of
Proposition~\ref{prop:forcing}
and in Lemma~\ref{mainerrorbound}, which are also $\le\delta_n$. 


Define 
\begin{equation}
x = \lim_{n\to\infty} x^{(n)}\in C^0(\zeroone,\ltwos);
\end{equation}
the limit exists  because the sequence $(x^{(n)})_{n\in\naturals}$
is constructed so as to be Cauchy in 
$C^0(\zeroone,\ltwos)$, as stated in \eqref{itiscauchy}.
Together, \eqref{initialbigness} and \eqref{itiscauchy} ensure that
the component $x_0$ does not vanish identically as a function of 
$t\in\zeroone$, so 
$x$ is a nonzero element of $C^0(\ltwos)$.
Because $x^{(n)}(0)\equiv 0$, the same holds for $x$;
that is, $x$ satisfies the desired initial condition at time $t=0$.

\subsection{Existence of $\scriptn(x)$}
In order to show that $x$ satisfies the desired differential equation,
we must first show that $\scriptn(x)$ is well-defined.
\begin{lemma} \label{nonlineartermexists}
Let $(\frakm_N)$ be a uniformly bounded sequence of finitely supported functions
from $\integers$ to $\complex$, and suppose that
$\lim_{N\to\infty}\frakm_N(n)=1$ for all $n\in\integers$.
Define\footnote{
An abuse of notation; this $\scriptp_N$ 
is related to what is called
$\scriptp_N$ elsewhere in the paper by conjugation with the spatial
Fourier transform.}
the operators  $\scriptp_N y_n = \frakm_N(n)y_n$.
Then
$\lim_{N\to\infty} \scriptn(\scriptp_N x)$ 
exists in $C^{-1}(\zeroone,\ltwos)$ norm.
\end{lemma}

Two facts will be repeatedly used in the proof of 
Lemma~\ref{nonlineartermexists}. Firstly,
\begin{equation} \label{scriptncontinuity}
\|\scriptn(v)-\scriptn(w)\|_{C^{0}(\zeroone,\ell^1)}
\le C
\|v-w\|_{C^{0}(\zeroone,\ell^1)}
\cdot \big(\|v\|_{C^{0}(\zeroone,\ell^1)}
+
\|w\|_{C^{0}(\zeroone,\ell^1)}\big)^2.
\end{equation}
Secondly, the operators $\scriptp_N$ are uniformly bounded
on $C^r(\zeroone,\scripth)$ for $\scripth=\ell^1$ and
$\scripth=\ltwos$, for $r=0$ and $r=-1$.

For any $N,k$,
\begin{equation}
\|\scriptn(\scriptp_N x^{(k)})-\scriptn(\scriptp_N x)\|_{C^0(\zeroone,\ell^1)}
\le CN^3 2^{-k}, 
\end{equation}
since $x^{(k)}-x\le 2^{-k}$ in $C^0(\zeroone,\ell^\infty)$ norm
and $\scriptp_N y$ is supported $[-3N,3N]$ for any $y$.
Thus for any index $J$ 
\begin{equation} \label{telescope}
\scriptn(\scriptp_N x) = 
\scriptn(\scriptp_N x^{(J)})
+ \sum_{j=J}^\infty 
\big[
\scriptn(\scriptp_N x^{(j+1)}) - \scriptn(\scriptp_N x^{(j)})
\big]
\end{equation}
with convergence in the $C^0(\zeroone,\ell^1)$ norm.

For any fixed $k$,
$\scriptp_Nx^{(k)}\to x^{(k)}$ in $C^0(\zeroone,\ell^1)$
norm since the multipliers $\frakm_N$ are uniformly bounded
and tend pointwise to $1$.
Therefore 
\begin{equation} \label{easyconvergence}
\|
\scriptn(\scriptp_N x^{(k)}) - \scriptn(x^{(k)})
\|_{C^{0}(\zeroone,\ell^1)}
\to 0
\text{ as $N\to\infty$}
\end{equation}
by \eqref{scriptncontinuity}.

\begin{lemma} \label{nuisancelemma}
If the construction of the sequence $(x^{(n)})_{n\in\integers}$
is carried out so that each $\delta_n$ is sufficiently
small relative to quantities determined at earlier steps of the construction,
then there exists $C<\infty$ such that for all $k$ and all $N$,
\begin{equation}
\|\scriptn(\scriptp_N x^{(k+1)}) 
- \scriptn(\scriptp_N x^{(k)})\|_{C^{-1}(\zeroone,\ltwos)}
\le C2^{-k}.
\end{equation}
\end{lemma}

\begin{proof}
Rewrite 
\begin{multline} \label{rewritten}
\scriptn(\scriptp_N x^{(k)} + \scriptp_N h^{(k)}) 
-\scriptn(\scriptp_N x^{(k)})
\\
=
\big[\scriptn(\scriptp_N x^{(k)} + \scriptp_N h^{(k)})
-\scriptn(\scriptp_N x^{(k)}) - \scriptn(\scriptp_N h^{(k)}) \big]
+ \nmain(\scriptp_N h^{(k)})
+ \ndiag(\scriptp_N h^{(k)}).
\end{multline}
Now $\ndiag(\scriptp_N h^{(k)})$ 
can be bounded in $C^0(\ltwos)$ norm exactly
as was done for $\ndiag(h^{(k)})$ in the proof of
Lemma~\ref{mainerrorbound}, up to an additional factor of
$\|\frakm_N\|_{\ell^\infty}^3$.
The same applies to $\scriptn(\scriptp_N x^{(k)} + \scriptp_N h^{(k)})
-\scriptn(\scriptp_N x^{(k)}) - \scriptn(\scriptp_N h^{(k)})$
in comparison with $\scriptn( x^{(k)} +h^{(k)})
-\scriptn(x^{(k)}) - \scriptn(h^{(k)})$.

$\nmain(\scriptp_N h^{(k)})$
breaks up into two parts. First there is the contribution 
of all $3$-tuples $(k,l,m)\in S^\dagger{}^3$, with $l\ne k,m$, 
that are not of the form $(m_j,m'_j,m_j)$. The same analysis
given for $\nmain(h^{(k)})-f^{(k)}$ in the proof of 
Lemma~\ref{mainerrorbound}
applies to the sum of these terms, up to the factor of
$\|\frakm_N\|_{\ell^\infty}^3$. 
Thus the sum of these terms
is again as small as desired in $C^0(\ltwos)$ norm.

There remains the contribution of all $3$-tuples
$(m_j,m'_j,m_j)$. 
Any such $3$-tuple contributes exactly
$\frakm_N(m_j)^2\overline{\frakm_N(m'_j)}$ times
$f_{n_j}^{(k)}(t)$. Because $\frakm_N(m_j)^2\overline{\frakm_N(m'_j)}$
is independent of $t$, we therefore have the same upper bound
in $C^{-1}(\ltwos)$ as for $f^{(k)}$ itself, up to the factor 
$\|\frakm_N\|_{\ell^\infty}^3$. 
\end{proof}

Lemma~\ref{nonlineartermexists} follows directly
from the combination of Lemma~\ref{nuisancelemma}
with \eqref{telescope} and \eqref{easyconvergence}.
Henceforth $\scriptn(x)$ is well-defined, via Lemma~\ref{nonlineartermexists}.

\begin{corollary}
\begin{equation}
\scriptn(x^{(k)})\to\scriptn(x) 
\text{ in $C^{-1}(\zeroone,\ltwos)$ norm as $k\to\infty$.}
\end{equation}
\end{corollary}

\begin{proof}
By Lemma~\ref{nuisancelemma},
\begin{equation}
\|\scriptn(\scriptp_N x)-\scriptn(\scriptp_N x^{(k)})
\|_{C^{-1}(\zeroone,\ltwos)}
\le C2^{-k},
\end{equation}
uniformly in $N$.
By Lemma~\ref{nonlineartermexists} and \eqref{easyconvergence},
\begin{equation}
\scriptn(\scriptp_N x)-\scriptn(\scriptp_N x^{(k)})
\to
\scriptn(x)-\scriptn( x^{(k)})
\text{ in ${C^{-1}(\zeroone,\ltwos)}$ norm as $N\to\infty$.}
\end{equation}
Therefore
\begin{equation}
\| \scriptn(x)-\scriptn( x^{(k)})\|_{C^{-1}(\zeroone,\ltwos)}
\le C2^{-k}.
\end{equation}
\end{proof}

\subsection{A solution of the Cauchy problem}
By definition of $f^{(n)}$,
$x^{(n)}(t) = \int_0^t \scriptn(x^{(n)}(s))\,ds + \int_0^t f^{(n)}(s)\,ds$.
Since $f^{(n)}\to 0$ in $C^{-1}(\zeroone,\ltwos)$ norm,
$x^{(n)}\to x$ in $C^0(\zeroone,\ltwos)$,
and
$\scriptn(x^{(n)})\to \scriptn(x)$ in $C^{-1}(\zeroone,\ltwos)$,
it follows at once that
\begin{equation} \label{ODEforx}
x(t) = \int_0^t \scriptn(x(s))\,ds.
\end{equation}

Define $u\in C^0(\zeroone,H^s)$ by 
\begin{equation}\widehat{u}(t,n) = e^{-in^2 t}x_n(t).\end{equation}
Since $x^{(n)}(t,x)\in C^0(\zeroone,\ltwos)$ vanishes identically
for $t=0$ and tends to $x$ in $C^0(\zeroone,\ltwos)$ norm,
$u$ satisfies the initial condition $u(0,\cdot)\equiv 0$.
Lemma~\ref{nonlineartermexists} states 
in equivalent form that $\scriptn(u)$
exists in the sense of Definition~\ref{defn:truncation}.
\eqref{ODEforx} implies that
$u$ is a weak solution in the extended sense
of the modified nonlinear Schr\"odinger equation.
This completes the proof of Theorem~\ref{thm:first}.
\qed

\section{Variants}
\subsection{The analogue for $\scripth^p$}

The proof of Theorem~\ref{thm:second} 
is quite similar to that of Theorem~\ref{thm:first}. The only 
significant change arises
in the proof of Proposition~\ref{prop:forcing},
for one cannot make $\|h\|_{C^0(\zeroone,\ell^p)}$
arbitrarily small simply by selecting $S^\dagger\subset[M,\infty)$
for $M$ arbitrarily large, 
as can be done for $\|h\|_{C^0(\zeroone,\ltwos)}$. 

The key now is that
with a modification of the set
$S^\dagger$ of spatial Fourier modes in the support of the new
driving force $g$, 
making $\nmain(h)\approx f$ 
requires a lower bound on $h$ in $C^0(\zeroone,\ell^2)$
but not in $C^0(\zeroone,\ell^p)$ for $p>2$.
Let $S=\{n_j: 1\le j\le A\}$ be as in 
the proof of Proposition~\ref{prop:forcing}.
$S^\dagger$ will now be taken to consist of
elements $m_0$ and $m_{j,i},m'_{j,i}$ for $1\le j\le A$
and $1\le i\le K$ where the new parameter $K$ is to be determined.
A large integer $m_0\in\naturals$ is chosen first,
then $m_{1,1}<m_{1,2}<\cdots<m_{1,K}<m_{2,1}<m_{2,2}<\cdots<m_{2,K}<m_{3,1}
<\cdots<m_{A,K}$ are chosen in that order, each sufficiently
large relative to all its predecessors for later purposes,
and then the quantities
$m'_{j,i}$ are uniquely determined by the relations
\begin{equation}
m_0+ m_{j,i}- m'_{j,i}=n_j
\text{ for all } j,i.
\end{equation}
If $m_0$ is chosen so that $m_0-n_j>0$ for all $j$ then
there is no obstruction to choosing $m_{j,i},m'_{j,i}$
so that this equation holds and $m_0,m_{j,i},m'_{j,i}$ are three
distinct integers.

$h_{m_0}(t)$ is defined to be the constant function $c\eps$
where $\eps$ is the small quantity in the conclusion of
the Proposition, and $c$ is some sufficiently small fixed constant.
Coefficients $\{h_m: m\in S^\dagger\}$ are chosen to be
$C^1$ functions satisfying
\begin{gather}
i\omega h_{m_{j,i}}\overline{h_{m'_{j,i}}}h_{m_0}(t)
\equiv \tfrac12 K^{-1} e^{-i\sigma(m_j,m'_j,m_j,n_j)t}
f_{n_j}(t)
\\
\|h_m\|_{C^0(\zeroone)}
\le C\big(\eps^{-1} K^{-1}\|f_{n_j}\|_{C^0} \big)^{1/2}
\ \text{ if $m=m_{j,i}$ or $m=m'_{j,i}$.}
\label{Kgain}
\end{gather}
for each $1\le j\le A$ and each $1\le i\le K$.
If $p$ is strictly greater than $2$ then for any given $\delta>0$,
$\{h_{j,i}\}$ can be made to satisfy
\begin{equation}
\big(\sum_{m\ne m_0\in S^\dagger} \|h_{m}\|_{C^0}\big)^{1/p}
\le\delta,
\end{equation}
by choosing $K$ to be sufficiently large 
as a function of $\eps,\delta$,
for the factor of $K^{1/p}$ arising from the number of terms on
the left-hand side is more than compensated for by the factor of $K^{-1/2}$
in \eqref{Kgain},
and this allows us to absorb the factor $\eps^{-1/2}$ in \eqref{Kgain}.
The remainder of
the proof of Proposition~\ref{prop:forcing} is unchanged.
Repeated applications of the Proposition establish Theorem~\ref{thm:second},
just as for Theorem~\ref{thm:first}.
\qed

\subsection{Quadratic nonlinearities}
Consider the nonlinearity $Q(u) = u^2$; the 
discussion will apply to $\bar u^2$ and $|u|^2-\mu(|u|^2)$ 
with very minor changes which are left to the reader.
If $S=\{n_j: 1\le j\le A\}$
then we set $S^\dagger=\{m_j,m'_j: 1\le j\le A\}$
where $m_j+m'_j=n_j$ and
$|m_j|,|m'_j|\ge M$ for all $j$.
The conditions on $\{h_m\}$ now become
\begin{equation} \label{quadraticbeqn}
i\omega h_{m_j}(t)h_{m'_j}(t)
\equiv\tfrac12 e^{-i(m_j^2+{m'_j}^2-n_j^2)t}f_{n_j}(t).
\end{equation}
By choosing $m_1$ sufficiently large and then
$|m_j|$ sufficiently large relative to $|m_{j-1}|$  
we may ensure that the analogue of Lemma~\ref{mainerrorbound} holds.
The rest of the argument is unchanged.
\qed



\begin{thebibliography}{20}
\bibitem{biagioni}
H.~A.~Biagioni,
{\em A nonlinear theory of generalized functions},
Second edition. Lecture Notes in Mathematics, 1421.
Springer-Verlag, Berlin, 1990.

\bibitem{bourgainperiodic}
J.~Bourgain,
{\em Fourier transform restriction phenomena for certain lattice subsets
and applications to nonlinear evolution equations. I. Schr\"odinger equations},
Geom.\ Funct.\ Anal.\  3  (1993),  no.\ 2, 107--156.
MR1209299 (95d:35160a)

\bibitem{bourgainkdv}
\bysame,
{\em Fourier transform restriction phenomena for certain lattice subsets 
and applications to nonlinear evolution equations. II. The KdV-equation},  
Geom. Funct. Anal.  3  (1993),  no. 3, 209--262.
MR1215780 (95d:35160b) 

\bibitem{burqgerardtzvetkov}
N. Burq, P. G\'erad\ and\ N. Tzvetkov, An instability property of the nonlinear Schr\"odinger equation on $S\sp d$, Math. Res. Lett. {\bf 9} (2002), no.~2-3, 323--335. MR1909648 (2003c:35144)

\bibitem{christpowerseries}
M.~Christ,
{\em Power series solution of a nonlinear Schr\"odinger equation},
preprint December 2004.

\bibitem{cct1}
M.~Christ, J.~Colliander, and T.~Tao,
{\em Asymptotics, frequency modulation, and low regularity ill-posedness 
for canonical defocusing equations},
Amer. J. Math. 125 (2003), no. 6, 1235--1293.  MR2018661 (2005d:35223) 

\bibitem{cct2}
\bysame
{\em Illposedness for nonlinear Schr\"odinger and wave equations},
to appear, Annales IHP Analyse Non Lin\'eaire.

\bibitem{cct3}
\bysame,
{\em Instability of the periodic nonlinear Schr\"odinger equation}, preprint,
math.AP/0311227.


\bibitem{colombeau}
J.-F.~Colombeau,
{\em Multiplication of distributions.
A tool in mathematics, numerical engineering and theoretical physics},
Lecture Notes in Mathematics, 1532. Springer-Verlag, Berlin, 1992.


\bibitem{dix}
D.~B.~Dix,
{\em Nonuniqueness and uniqueness in the initial-value problem for Burgers' equation}, 
SIAM J. Math. Anal.  27  (1996), no. 3, 708--724. MR1382829 (97c:35174)

\bibitem{planchonetal}
G.~Furioli, F.~Planchon, and E.~Terraneo, 
{\em Unconditional well-posedness for semilinear Schr\"odinger 
and wave equations in $H\sp s$},  
Harmonic analysis at Mount Holyoke (South Hadley, MA, 2001),  
147--156, Contemp. Math., 320, Amer. Math. Soc., Providence, RI, 2003,
MR1979937.

\bibitem{kappelertopalovkdv}
T.~Kappeler and P.~Topalov,
{\em Global Well-Posedness of KdV in $H^{-1}(T,R)$}, preprint.

\bibitem{kenigetal}
C.~E.~Kenig, G.~Ponce, and L.~Vega,
{\em On the ill-posedness of some canonical dispersive equations},  
Duke Math. J.  106  (2001),  no. 3, 617--633.

\bibitem{kpvquad}
\bysame,
{\em Quadratic forms for the $1$-D semilinear Schr\"odinger equation}, 
Trans. Amer. Math. Soc.  348  (1996), no. 8, 3323--3353. MR1357398 (96j:35233)


\bibitem{scheffer}
V.~Scheffer,
{\em 
An inviscid flow with compact support in space-time},  
J. Geom. Anal. 3 (1993), no. 4, 343--401.

\bibitem{shnirelman}
A.~Shnirelman,
{\em On the nonuniqueness of weak solution of the Euler equation},
Comm. Pure Appl. Math. 50 (1997), no. 12, 1261--1286.

\end{thebibliography}
\end{document}